\documentclass[11pt]{amsart}
\usepackage{amssymb}
\usepackage{amsfonts}
\usepackage{amscd}
\usepackage{amsmath}
\usepackage{mathrsfs}
\usepackage{curves}
\usepackage{epsfig}
\usepackage[pass]{geometry}
\usepackage{graphicx}
\usepackage{verbatim}
\usepackage{latexsym}
\usepackage{eucal}
\usepackage{enumerate} 
\usepackage{hyperref}
\usepackage{xcolor}
\usepackage{caption}
\usepackage{mathtools}

\newtheorem{theorem}{Theorem}[section]

\def \mcb {{\mathcal B}}

\def \mcd {{\mathcal D}}

\def \mck {{\mathcal K}}

\def \mcs {{\mathcal S}}

\def \mcv {{\mathscr V}}

\def \mbn {{\mathbb N}}
\def \mbr {{\mathbb R}}
\def \mbs {{\mathbb S}}

\def \beqq {\begin{equation}}
\def \eeqq {\end{equation}}
\def \WF {\text{WF}}

\def \bpf {\begin{proof}}
\def \epf {\end{proof}}
\def \beq {\begin{equation*}}
\def \eeq {\end{equation*}}

\def \eps {\epsilon}   
   
\def \La {\Lambda}

\def \p {\partial}
\def \ha {\frac{1}{2}}

\def \tilde {\widetilde}

\begin{document}
\title[Robustness in limited angle tomography]{Robustness of data-driven approaches in limited angle tomography}
\author{Yiran Wang and Yimin Zhong}
\address{Yiran Wang
\newline
\indent Department of Mathematics, Emory University}
\email{yiran.wang@emory.edu}
\address{Yimin Zhong
\newline
\indent Department of Mathematics and Statistics, Auburn University}
\email{yzz0225@auburn.edu}
\begin{abstract} 
The limited angle Radon transform is notoriously difficult to invert due to its ill-posedness. In this work, we give a mathematical explanation that data-driven approaches can stably reconstruct more information compared to traditional methods like filtered backprojection. In addition, we use experiments based on the U-Net neural network to validate our theory.  %More precisely, we give a condition on the data set under which the data-driven approach could potentially reach robustness of Lipschitz type. 
\end{abstract}

\maketitle 

%%%%%%%%%%%%%%%%%%%
\section{Introduction}\label{sec-intro} 
Consider the limited angle Radon transform on $\mbr^2$. Let $(x_1, x_2)$ be the coordinate for $\mbr^2$. For $\delta \in (0, \pi/2),$ let  $\mcs_\delta = \{(\cos\theta, \sin \theta)\in \mbs^1:  |\theta| \leq \delta\}.$  We use the parallel beam parametrization of the Radon transform. Let $v\in \mcs_\delta$ and $v^\perp$ be the unit vector orthogonal to $v$ so that $\{v, v^\perp\}$ has the same orientation as the coordinate frame. The limited angle Radon transform of a function $f$ on $\mbr^2$ is 
\beqq\label{eq-xray}
R_\delta f(p, v) = \int_{\mbr} \chi_\delta(v) f(p v + tv^\perp)dt, \quad  v\in \mbs^1, p\in \mbr
\eeqq 
where $\chi_\delta$ is the characteristic function of $\mcs_\delta$ in $\mbs^1.$ For simplicity, we assume that $f$ is supported in the unit disk $\mcb  = \{x\in \mbr^2: |x|<1\}$ where $|\cdot|$ stands for the  Euclidean norm. It is a classical result that $R_\delta$ is injective on compactly supported distributions, for instance, \cite{Nat, StUh}.  However, 
 the inverse problem of recovering $f$ from $R_\delta f$ is severely ill-posed. This was demonstrated from the behavior of the singular values in \cite{Dav}. On the microlocal level, it is clear that some information of $f$ is lost. In particular, Quinto in~\cite{Qui} gave a microlocal description of the singularities of $f$ that can be recovered, called ``visible singularities", and those that cannot, called ``invisible singularities", see Section 2 for more details. Due to the ill-posedness, traditional methods such as filtered back-projection or iterative methods for solving a variational regularized problem (e.g., a sparsity enforcing approach solved with ISTA, see e.g.~\cite{BGR, Fri}) often produce artifacts in the reconstruction, and it is very challenging to recover the invisible information.

Recently, there have been increasing efforts in applying data-driven approaches such as deep neural networks to limited angle tomography. Among the already vast literature, we mention \cite{BAK, BKL, BGR, GMJ, Goy, GuYe, RMB} which are more relevant to our work, and \cite[Section II.B]{Pan} for more references. Multiple numerical experiments reported that data-driven approaches can significantly reduce the artifact phenomena and even learn the invisible singularities. The success seems difficult to understand from the existing mathematical theory. Also, the robustness of data-driven approaches is a concern, see for example \cite{Hua}. It is suspected that data-driven approaches should be tailored to specific tasks, not for general image reconstruction. But a precise understanding is still missing.

\begin{figure}[t]
\centering
\includegraphics[scale=0.55]{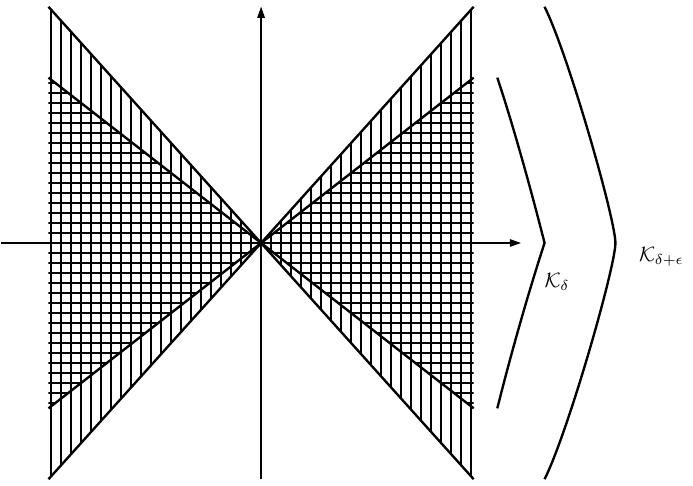}
\caption{Illustration of the increased stability in the Fourier space. Traditional methods can recover $\hat f$ in the cone $\mck_\delta$ (defined in Section \ref{sec-stab}). Data-driven approaches can recover $\hat f$ in a larger cone $\mck_{\delta + \eps}$. The size of the cone depends on the specific training set. }
\label{fig-fourier}
\end{figure}

The goal of this paper is to provide a mathematical explanation. It is useful to interpret the instability of the limited angle Radon transform from the following viewpoint: there is no stability estimate of the form 
\beqq\label{eq-stab0}
\|f\|_{H^s(\mbr^2)} \leq C(s, t) \|R_\delta f\|_{H^{s+t}(\mbr\times \mbs^1)}
\eeqq
for any $s, t\in \mbr$, see Section VI.2 of \cite{Nat}. This is because, roughly speaking, after taking the limited angle Radon transform, the Fourier transform  $\hat f$ of $f$ is lost in a conic region (corresponding to the missing angles), see Figure \ref{fig-fourier}. Recovery of $\hat f$ in the rest of the Fourier space ($\mck_\delta$ in Figure \ref{fig-fourier}) however is stable, see Theorem \ref{thm-main0}. This result is sharp. For data-driven approaches, we consider $f$ in a training data set $\mcd$ and show that $\hat f$ can be stably recovered in a larger cone ($\mck_{\delta +\eps}$ in Figure \ref{fig-fourier}), see Theorem \ref{thm-main}. The size of the cone depends on $\mcd$ and it is possible to cover the whole Fourier space for a certain $\mcd$. This explains the phenomena of ``learning the invisible". We remark that our theory does not rely on the specific structure of the neural network, at least for supervised learning. In Section \ref{sec-num}, we use a simple U-Net to validate our theory.  The concluding summary is offered in Section \ref{sec-conclude}.

%%%%%%%%%%%%%%%%%%%
\section{The visible and invisible singularities}\label{sec-sing} 
In \cite{Qui1}, Quinto described some principles for recovering features of an object from limited X-ray CT. Roughly speaking, if a boundary of a feature of the body is tangent to a line in a limited data set, then that boundary should be easy to reconstruct (called visible boundaries). Otherwise, the boundary should be difficult to reconstruct from the limited data (i.e. invisible boundaries).  We use the FBP reconstruction to demonstrate the principle. 

Let $\theta\in (-\delta, \delta)$ and define the limited angle back-projection for $R_\delta$ as
\beq
R^*_\delta g(x) = \int_{-\delta}^{\delta} g(v(\theta), x\cdot v(\theta)) d\theta. 
\eeq
Let $\La$ be the Riesz potential 
\beq
\La g(v, p) = \frac{1}{2\pi} \int_{-\infty}^\infty \int_{-\infty}^\infty e^{i\tau (p - s)} |\tau| g(v, s) ds d\tau.
\eeq
Then we consider the following FBP for limited angle transform
\beqq\label{eq-fbp}
L_\delta f = \frac{1}{4\pi} R_\delta^* \La R_\delta f.
\eeqq
Note that this is not an inversion formula. We describe what can be reconstructed from this formula. 

We introduce some notions from microlocal analysis. Let $f$ be a function on $\mbr^2$ and let $x_0\in \mbr^2, \xi_0\in \mbr^2\backslash \{0\}$. By definition, $f$ is smooth at $x_0$ in direction $\xi_0$ if there is a smooth function $\phi\in C_c^\infty(\mbr^2)$ with $\phi(x_0)\neq0$ and an open cone $\mcv$ containing $\xi_0$ such that for any $N\in \mbn$ there is $C_N$ such that 
\beq
|\widehat{\phi f}(\xi)| \leq C_N (1 + |\xi|)^{-N}, \text{ for all } \xi \in \mcv,
\eeq
where $\hat{}$ denotes the Fourier transform in $\mbr^2.$ 
If $f$ is not smooth at $x_0$ in direction $\xi_0$, then $(x_0, \xi_0)\in \WF(f)$ the wavefront set of $f$. Let $(x_0, \xi_0)\in \WF(f).$ We say $(x_0, \xi_0)$ is a {\em visible singularity} of $f$ (for $L_\delta$) if $\xi_0$ is parallel to some $v\in \mcs_\delta$. We say $(x_0, \xi_0)$ is an {\em invisible singularity} of $f$ if $\xi_0$ is not parallel for any $v\in \mcs_\delta$. We have 
\begin{theorem}[Theorem 2 of \cite{Qui1}]\label{thm-qui1}
Let $f$ be a function of compact support and $(x_0, \xi_0)\in \WF(f)$.
\begin{enumerate}
\item If $(x_0, \xi_0)$ is a visible singularity of $f$ then $(x_0, \xi_0)\in \WF(L_\delta f)$.
\item  If $(x_0, \xi_0)$ is an invisible singularity of $f$ then $(x_0, \xi_0)\notin \WF(L_\delta f)$.
\end{enumerate}
\end{theorem}

We remark that it is possible to obtain the more precise microlocal Sobolev regularity of the visible singularities, see Theorem 3.1 of \cite{Qui1}. 
Within the microlocal framework, one can also explain the artifact phenomena. For our setup, the artifact can occur on lines parallel to $v\in \p\mcs_\delta$ and the lines are tangent to the boundary of some feature in the object, see Principle 2 of \cite{Qui}. %The artifact can be alleviated by introducing smooth cut-offs, see \cite{Qui1} for details. 
 
%%%%%%%%%%%%%%%%%%%
\section{The stable reconstruction}\label{sec-stab} 
We establish a quantitative estimate for recovering visible singularities. Our starting point is the simple injectivity proof for $R_\delta$, see \cite[Theorem 3.4]{Nat}. 
The transform $R_\delta f(p, v)$ is a function defined on $\mbr\times \mbs^1.$ We take the Fourier transform of $R_\delta f$ in $p$ to get  
\beqq\label{eq-fourier}
\begin{aligned}
\widehat{R_\delta f}(\tau, v) &= \int_\mbr \int_\mbr e^{-\imath p\tau} \chi_\delta(v) f(p v + tv^\perp) dt dp \\
 &= \int_\mbr \int_\mbr e^{-\imath (pv  + tv^\perp) \cdot \tau v } \chi_\delta(v) f(p v  + tv^\perp) dt dp  =   \hat{f}(\tau v) \chi_\delta(v).
\end{aligned}
\eeqq
Here, $\imath^2 = -1$. Because $R_\delta f(p, v) = 0$, we get $\hat{f}(\xi)  = 0$ for $\xi\in \mbr^2$ parallel to $v$. Thus $\hat f(\xi) = 0$ in a cone  $\mck_\delta = \{\xi\in \mbr^2: \xi/|\xi|\in \mcs_\delta\}$ with non-empty interior for $\delta >0$. Because $\hat{f}$ is analytic in $\xi$, we conclude that $\hat{f}(\xi) = 0$ for all $\xi\in \mbr^2.$ So $f = 0$ after taking the inverse Fourier transform.  

We can further derive some stability estimates. Let $\tilde \chi_\delta$ be the characteristic function of $\mck_\delta$ in $\mbr^2$. Note that $\tilde \chi_\delta(\xi) = \chi_\delta(\xi/|\xi|), \xi\neq 0$. For simplicity, we use the same notation $\chi_\delta$ for both characteristic functions.  It should be clear from the context which one we are referring to. Also, we use $\chi_\delta(D)$ for the Fourier multiplier of $\chi_\delta(\xi)$ defined by $\chi_\delta(D)f =  (\chi_\delta \hat f)^\vee$ where ${}^\vee$ denotes the inverse Fourier transform.  

\begin{theorem}\label{thm-main0}
For $\delta \in (0, \pi/2]$ and $f\in H^{-1/2}(\mbr^2),$ we have
\beqq\label{eq-mstab}
\|\chi_\delta(D)f\|_{H^{-1/2}(\mbr^2)}\leq   \|R_\delta f\|_{L^2(\mbr\times \mbs^1)}  .
\eeqq 
\end{theorem}
\bpf
Applying Plancherel's theorem and using \eqref{eq-fourier}, we get 
\beqq\label{eq-fourier2}
\begin{aligned}
 \| R_\delta f\|^2_{L^2(\mbr\times \mbs^1)} =& \|\widehat{R_\delta f}\|^2_{L^2(\mbr\times \mbs^1)}  
 =   \int_{\mbs^1} \int_{\mbr}\chi_\delta(v) |\hat{f}(\tau v ) |^2 d\tau dv \\
 =& \int_{\mck_\delta} |\xi|^{-1} |\hat{f}(\xi)|^2 d\xi \geq  \|\chi_\delta(D)f\|_{H^{-1/2}(\mbr^2)}^2  .
\end{aligned}
\eeqq  
This finishes the proof.  
\epf 
This theorem is a quantitative version of Theorem \ref{thm-qui1} on stably recovering visible singularities of $f$. It is also clear that the stability estimate \eqref{eq-stab0} fails because of $(1 - \chi_{\delta}(D))f$ which sits in the kernel of $R_\delta.$ 

We seek improvement of the estimate to include invisible singularities. For $\eps>0$ small we define $\chi_{\delta+\eps}(D)$ as a Fourier multiplier of $\chi_{\delta+\eps}(\xi), \xi\in \mbr^2$. For $N\geq 0$, we define
\beqq\label{eq-dn}
\begin{gathered} 
\mcd_{N, \eps} = \{f\in L^2(\mcb):  \|f\|_{L^2} \leq N \|\chi_{\delta + \eps}(D)  f\|_{H^{-1/2}} \}.
\end{gathered}
\eeqq
Here, we recall that $\mcb$ is the unit disk in $\mbr^2.$ Observe that $L^2(\mcb) = \cup_{N\geq 0}\mcd_{N, \eps}.$ 
The following theorem says that for functions in $\mcd_{N, \eps}$, the stability estimate \eqref{eq-mstab} can be improved. 

 \begin{theorem}\label{thm-main}
 For $N\geq 0$,  there is $\eps>0$ small such that 
\beqq\label{eq-stab1} 
\|\chi_{\delta + \eps}(D)f\|_{H^{-1/2}(\mbr^2)} \leq  2 \|R_\delta  f \|_{L^2(\mbr\times \mbs^1)} 
\eeqq 
for all $f\in \mcd_{N, \eps}.$
  \end{theorem} 
 
\bpf 
We decompose
\beqq\label{eq-decomp}
R_\delta f = R_\delta \chi_{\delta + \eps}(D)f + R_\delta (1 - \chi_{\delta +\eps}(D))f.   
\eeqq
Note that the second term on the right vanishes because $(1 - \chi_{\delta +\eps}(D))f$ is in the kernel of $R_\delta.$ So we estimate the first term. 
Using \eqref{eq-fourier}, we have  
\beqq\label{eq-deri0}
\begin{aligned} 
& R_{\delta}\chi_{\delta + \eps}(D)f(p, v) - R_{\delta+\eps} \chi_{\delta + \eps}(D)f(p, v)  \\
= &  (2\pi)^{-1} \int_\mbr e^{  \imath p \tau }(\chi_\delta(v)  - \chi_{\delta +\eps}(v)) \chi_{\delta + \eps}(\tau v)\hat f(\tau v) d\tau.
 \end{aligned}
\eeqq  
Note that $|\chi_\delta(v) - \chi_{\delta +  \eps}(v)|$ is bounded by one and supported in a set of $\mbs^1$ with measure $2\eps$. Also, the left-hand side of \eqref{eq-deri0} is compactly supported in $p$. Denote $U_{\lambda} = \{|\tau|\le \lambda\}\subset \mbr$, where we choose $\lambda = (2\pi\eps)^{-\ha}$, then using Plancherel's theorem again, we get  
\beqq\label{eq-deri1b}
\begin{aligned} 
& \| R_{\delta}\chi_{\delta + \eps}(D)f - R_{\delta+\eps} \chi_{\delta + \eps}(D)f  \|_{L^2}^2 \\
= & \int_{\mbr} \int_{\mbs^1} |(\chi_\delta(v)  - \chi_{\delta +\eps}(v)) \chi_{\delta + \eps}(\tau v)\hat f(\tau v)|^2  dv  d\tau \\ 
\leq & \int_{U_{\lambda}} \int_{\mbs^1} |(\chi_\delta(v)  - \chi_{\delta +\eps}(v))\hat f(\tau v)|^2  dv  d\tau + \frac{1}{\lambda}\|f\|_{L^2}^2\\ 
  \leq &  \pi \|f\|_{L^2}^2 \int_{U_{\lambda}}   \int_{\mbs^1}  (\chi_\delta(v)  - \chi_{\delta +\eps}(v))^2  dv   d\tau   + \frac{1}{\lambda}\|f\|_{L^2}^2 \\ 
\leq & (4\pi \eps \lambda   + \frac{1}{\lambda})\|f\|_{L^2}^2 
 \leq 4N^2 \sqrt{\pi \eps}  \|\chi_{\delta + \eps}(D)f\|^2_{H^{-1/2}}.
\end{aligned}
\eeqq
In the derivation, we used the fact that $\|\hat h\|_{L^\infty}\leq \sqrt{\pi} \|h\|_{L^2}$ for $h$ compactly supported in $\mathcal{B}$ and the inequality
\beq 
\begin{aligned}
\int_{|\tau|>\lambda } \int_{\mbs^1}    |\hat f(\tau v)|^2   dv d\tau &\leq \frac{1}{\lambda }\int_{|\tau| > \lambda } \int_{\mbs^1}   |\tau| |\hat f(\tau v)|^2  dv  d\tau  \le\frac{1}{\lambda}\|f\|_{L^2}^2 . 
\end{aligned}
\eeq
We can apply Theorem \ref{thm-main0} to derive that 
\beq
\begin{aligned}
\|\chi_{\delta + \eps}(D)f\|_{H^{-1/2}} &\leq   \|R_{\delta+ \eps}  \chi_{\delta + \eps}(D)f\|_{L^2} \\
&\leq   \| R_\delta \chi_{\delta + \eps}(D)f\|_{L^2} + \|R_\delta \chi_{\delta + \eps}(D)f - R_{\delta+ \eps}\chi_{\delta + \eps}(D)f\|_{L^2}  \\
&\leq   \| R_\delta \chi_{\delta + \eps}(D)f\|_{L^2} + 2N (\pi\eps)^{\frac{1}{4}}  \|\chi_{\delta + \eps}(D) f\|_{H^{-1/2}} .
\end{aligned}
\eeq
For $\eps >0$ sufficiently small so that $2N (\pi\eps)^{\frac{1}{4}} < \ha$, we obtain 
\beqq\label{eq-est1}
\|\chi_{\delta + \eps}(D)f\|_{H^{-1/2}(\mbr^2)} \leq  2 \|R_\delta   \chi_{\delta + \eps}(D)f \|_{L^2(\mbr\times \mbs^1)}  ,
\eeqq 
which finishes the proof. 
\epf

We remark that the choices of the Sobolev norms in the condition in \eqref{eq-dn} are for simplicity only. The analysis can be adapted to other Sobolev spaces by slightly modifying the proofs of Theorem \ref{thm-main0} and \ref{thm-main}. For general data-driven methods, if there exists an appropriate pair of $(N, \eps)$ such that the whole image data set (including training and testing) $\mcd \subset \mcd_{N, \eps} \subset H^{-1/2}$, then for any limited angle Radon transform measurement $F^{\ast} = R_{\delta} f^{\ast} + \eta$, $f^{\ast}\in \mcd$ and $\eta$ is the additive noise, the optimization 
\begin{equation}
    \min_{{f}\in  \mcd_{N,\eps}} \|R_{\delta}  f - F^{\ast}\|_{L^2}
\end{equation}
provides a stable reconstruction for $\chi_{\delta+\eps} (D) f^{\ast}$ in $H^{-1/2}(\mbr^2)$, which reveals the information in a wider cone $\mathcal{K}_{\delta+\eps}\supset \mathcal{K}_{\delta}$. Regarding the neural network approaches, there are extensive experiments~\cite{BAK, BKL, BGR, Goy, JinKH} with different structures (e.g. U-Net~\cite{unet}, It-Net~\cite{GMJ}). Typically, the training data set is generated from a low dimensional manifold $\mathcal{M}\subset \mathcal{D}_{N, \eps}$ and the neural network $\mathcal{N}_{\boldsymbol{ \theta}}$ tries to learn the inversion operator by finding the optimal weights $\boldsymbol{ \theta} \in \mbr^m$ that minimizes the loss function 
\begin{equation}
    \min_{\boldsymbol{ \theta} \in \mbr^m} \frac{1}{n}\sum_{j=1}^n\|\mathcal{N}_{\boldsymbol{ \theta}}(\mathcal{A}^{\dagger} F_j) - f_j\|_{L^2}^2
    \label{eq:loss}
\end{equation}
over the training data set $\{f_j\}_{j=1}^n$ and the corresponding measurements $\{F_j\}_{j=1}^n$. The operator $\mathcal{A}^{\dagger}$ is a preprocessing operator (e.g. FBP or a learnable layer). Suppose the neural network is well trained on this data set so that $\mathcal{N}_\theta$ is a good approximation of $\mathcal{A}^\dagger R_\delta$. In principle, this would happen when the network architecture has the universal approximation property (see e.g. \cite{Yar} for convolutional neural networks, and \cite{Hor} for feedforward networks), and the property is achieved through the training. Then the learned inversion operator $\mathcal{N}_{\theta}$ would potentially satisfy the improved stability estimate \eqref{eq-stab1} which can explain the stable reconstruction of certain invisible singularities observed in the experimental studies~\cite{Hua}.

%%%%%%%%%%%%%%%%%%%
\section{The experiment}\label{sec-num}
%========================%
\subsection{Data set}
\begin{figure}[!htbp]
\centering
\vspace{-0.3cm}
\includegraphics[scale=0.7]{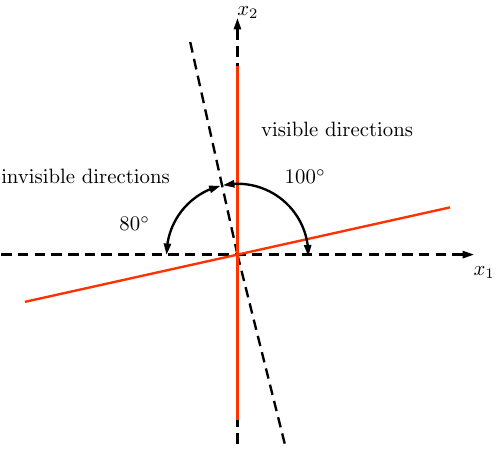}
\vspace{-0.3cm}
\caption{Illustration of the experimental setup and the corresponding visible/invisible directions according to the theory in Section \ref{sec-sing}. The red lines indicate the possible artifact directions.}
\label{fig-dir}
\end{figure}

\begin{figure}[!htb]
    \centering
    \includegraphics[width=0.42\textwidth]{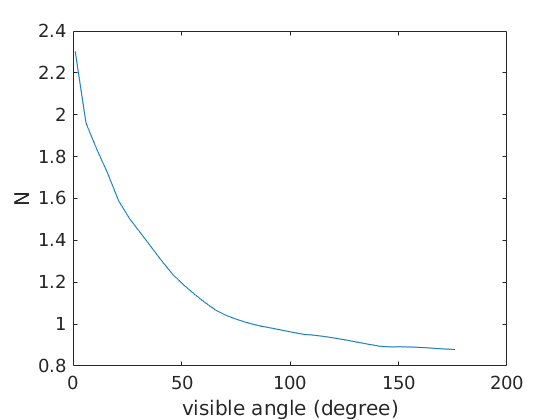}
    \includegraphics[width=0.42\textwidth]{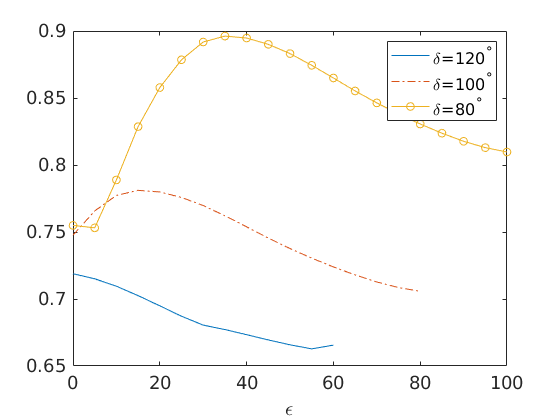}
    \caption{Left: Plot of $N_{\beta} = \sup_{f\in \mathcal{D}}\|f\|_{L^2}/\|\chi_{\beta}f\|_{H^{-1/2}}$ for the data set $\mathcal{D}$ with respect to each visible angle $\beta$. Right: Plots of $g(\eps)$ for visible angles $\delta=120^{\circ}, 100^{\circ}, 80^{\circ}$ (missing angles $\theta=60^{\circ}, 80^{\circ}, 100^{\circ}$). }
    \label{fig:N_estimate}
\end{figure}

We use the same data as in  \cite{BGR} which can be downloaded from \cite{data}. The data set consists of synthetic images of ellipses, where the number, locations, sizes, and intensity gradients of the ellipses are chosen randomly. The measurements are computed using the Matlab function \textsf{radon} with additive Gaussian noises. Also, the measurements are simulated at a higher resolution and then downsampled to an image of resolution $128\times 128$ to avoid the inverse crime. For the simulation below, we use measurements for missing angles $\theta = 60^\circ, 80^\circ, 100^\circ$ (out of $180^\circ$). See Figure \ref{fig-dir} for an example when $\theta=80^\circ$. 
The source code of the numerical experiments is hosted on GitHub\footnote{\href{https://github.com/lowrank/robustness-limited-angle-tomography}{https://github.com/lowrank/robustness-limited-angle-tomography}}. To estimate the parameters $N$ and $\eps$ for the data set $\mathcal{D}\subset \mathcal{D}_{N, \eps}$, we plot the curve for the critical value $N_{\beta} := \sup_{f\in \mathcal{D}}\|f\|_{L^2}/\|\chi_{\beta}f\|_{H^{-1/2}}$ in Figure~\ref{fig:N_estimate} after rescaling the images into a unit square. It shows the data set $\mathcal{D}\subset \mathcal{D}_{N, \eps}$ in~\eqref{eq-dn} for $N=1.2$ and all $\eps>0$ when the visible angle is greater than $60^{\circ}$. In addition, the curves of $g(\eps):=   \max_{f\in \mathcal{D}}\|\chi_{\delta+\eps} f\|_{H^{-1/2}} / \|R_{\delta} f\|_{L^2} $ with different $\delta$ are plotted in Figure~\ref{fig:N_estimate}, which implies the stability constant for the reconstruction should be relatively small. 

\subsection{The neural network}\label{subsec-net}
\begin{figure}[!hbt]
\centering
\includegraphics[scale=0.35]{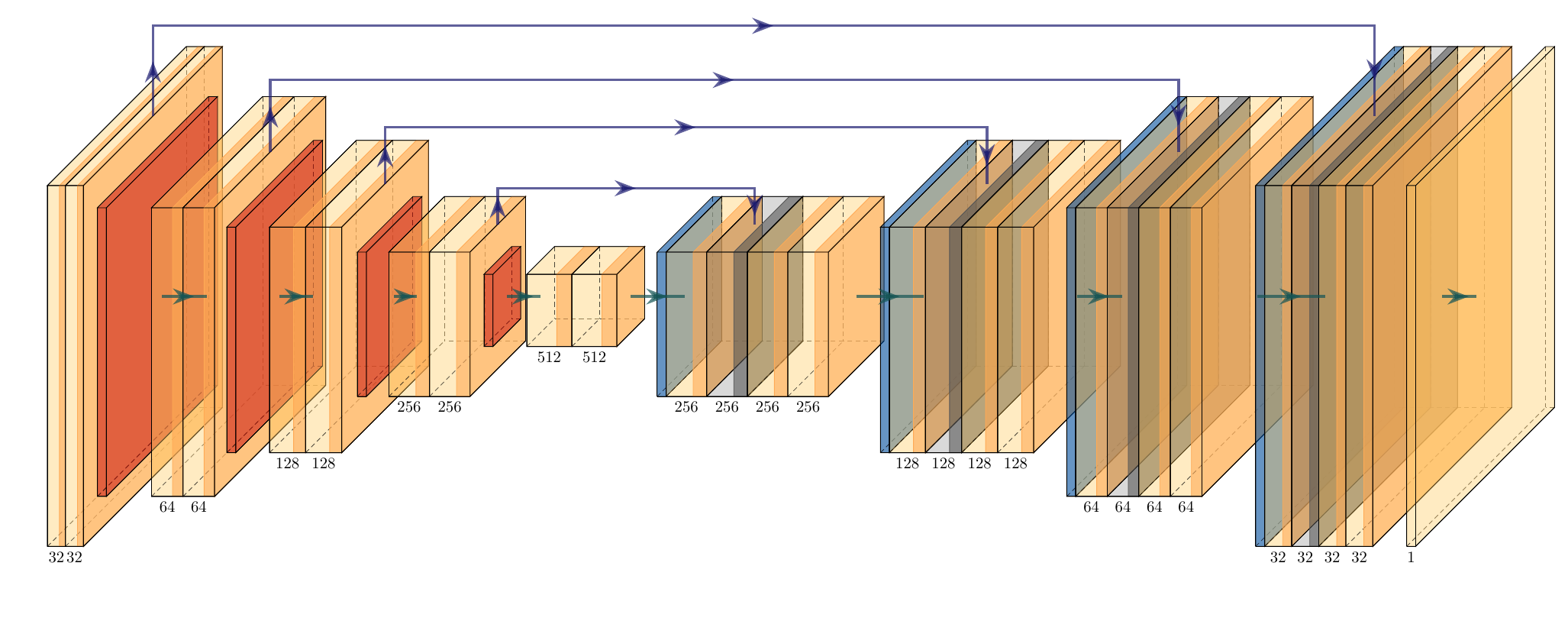}
\caption{Structure of the U-Net. The input is a $128\times 128$ image and the output is a $128\times 128$ image. The number of channels is indicated in the boxes. The arrows indicate the flow of the data.}
\label{fig: unet}
\end{figure}

\begin{figure}[!htb]
    \centering
    \includegraphics[width=0.325\linewidth]{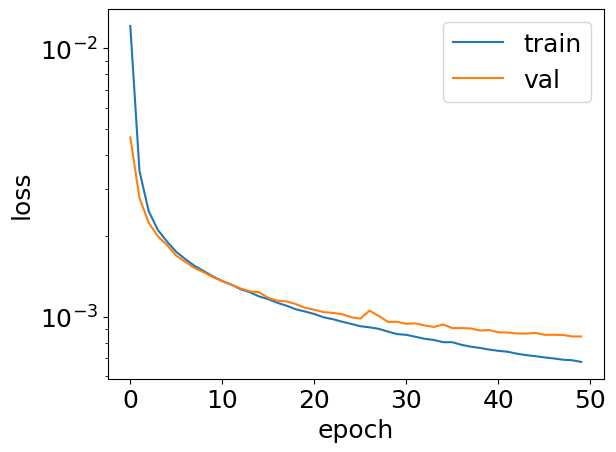}
    \includegraphics[width=0.325\linewidth]{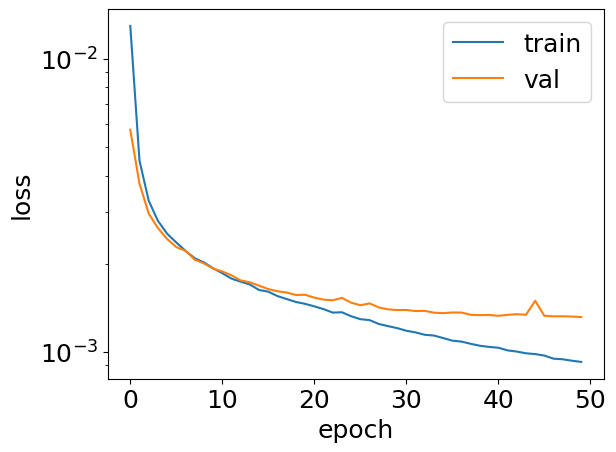}
    \includegraphics[width=0.325\linewidth]{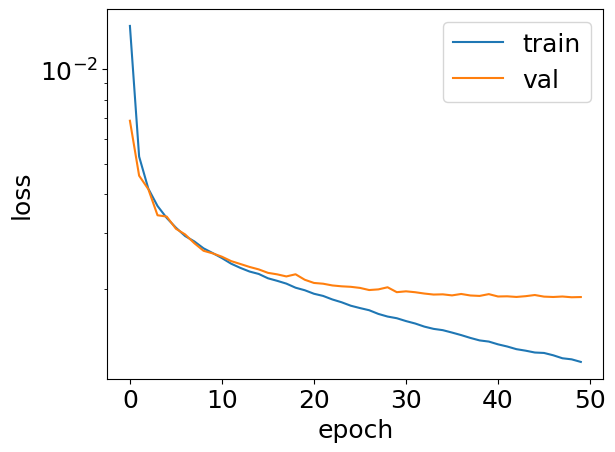}
    \caption{Loss~\eqref{eq:loss} of each epoch on training and validation sets. From left to right: missing angles $\theta=60^{\circ}, 80^{\circ}, 100^{\circ}$. }
    \label{fig-train_loss}
\end{figure}

\begin{figure}[!htbp]
\centering
\includegraphics[width=0.8\linewidth]{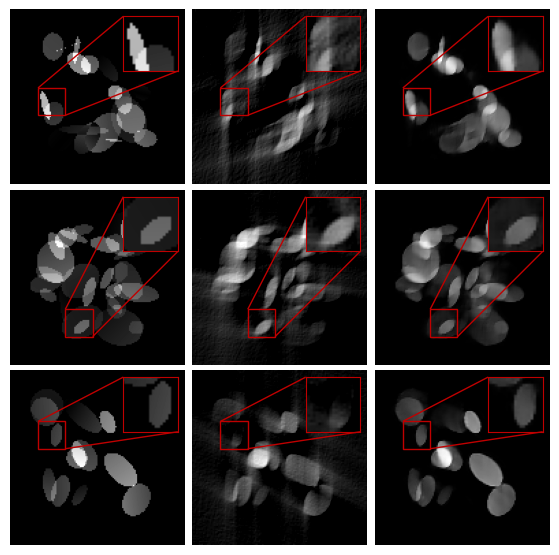}
% \includegraphics[scale=0.7]{figs_computed/fig_test522_0.jpeg} \quad 
% \includegraphics[scale=0.7]{figs_computed/fig_test522_1.jpeg} \quad
% \includegraphics[scale=0.7]{figs_computed/fig_test522_2.jpeg} 

% \vspace{0.1cm}
% \includegraphics[scale=0.7]{figs_computed/fig_test525_1.jpeg} \quad
% \includegraphics[scale=0.7]{figs_computed/fig_test525_2.jpeg}
\caption{Three numerical reconstructions using FBP and U-Net. Each row shows one example.  Each image's top-right shows a zoomed-in detail. From top to bottom: missing angle $\theta=100^{\circ}, 80^{\circ}, 60^{\circ}$. Left: the original image. Middle: the FBP reconstruction. Right: the U-Net reconstruction, from top to bottom, the ranges of pixel values for the reconstructions are $[-0.029, 1.036]$, $[-0.053, 0.950]$, and $[-0.047, 1.005]$. For display purposes, all image pixel values are clipped to $[0, 1]$. }
\label{fig-test}
\end{figure}

We take $\mathcal{A}^{\dagger}$ as the FBP inversion to preprocess the image data, which is fulfilled through the MATLAB function \textsf{iradon}. Then, we use a basic end-to-end U-Net~\cite{unet} for training with the network structure plotted in Figure~\ref{fig: unet}. 
The U-Net consists of three parts: downsampling,  bottleneck, and upsampling. The downsampling part has 4 blocks, connected with $2\times 2$ average pooling layers, and each block contains two convolutional layers with $3\times 3$ kernels and leaky ReLU activations. The upsampling part has 4 blocks, connected with upsampling layers, and each block contains two convolutional layers \emph{without} using activations. The final output convolution layer uses $1\times 1$ kernel. The channels of convolution layers are labeled in Figure~\ref{fig: unet}. We take $5000$ images as training samples among which $80\%$ are used for training and $20\%$ for validation, and the other $5000$ images as the test samples. The mini-batch size is $50$ and training uses the Adam optimizer with a small learning rate of $10^{-4}$ for 50 epochs. The training losses of the three cases (missing angle $\theta=60^{\circ}, 80^{\circ}, 100^{\circ}$) are plotted in Figure~\ref{fig-train_loss}. The discussions about the settings of the neural network can be found in Section~\ref{sec:dis}. 
We demonstrate three examples of the reconstruction on test images in Figure \ref{fig-test}. 
Observe that the examples' artifacts are almost completely removed, and some invisible singularities are reconstructed well. However, we still notice that the reconstructions are not perfect. Several thin ellipses with mainly invisible singularities are not reconstructed well in the examples.  
%========================%
\subsection{The robustness}\label{subsec-rob}
\begin{figure}[t]
\centering
\includegraphics[scale=0.7]{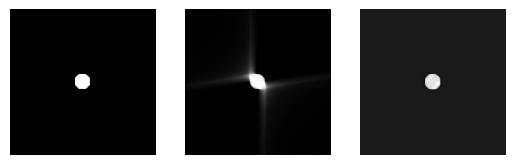}\\
\includegraphics[scale=0.7]{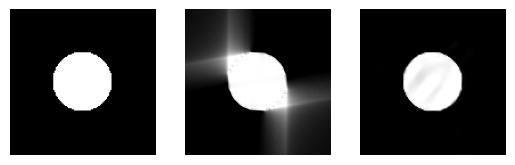}\\
\includegraphics[scale=0.7]{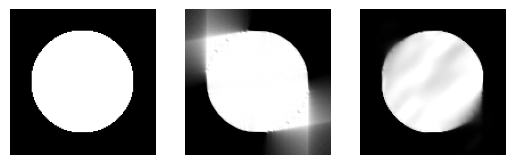}
\caption{Numerical test of the stability condition in \eqref{eq-dn} with missing angle $\theta=80^{\circ}$. Each row shows one example. Left: the original image. Middle: the FBP reconstruction. Right: the U-Net reconstruction, from top to bottom, the range of pixel values for the reconstructions are $[-0.12, 1.12]$, $[-0.08, 1.18]$, and $[-0.07, 1.17]$. For display purposes, all image pixel values are clipped to $[0, 1]$.}
\label{fig-disk}
\end{figure}

We consider the violation of the condition in \eqref{eq-dn}. For $a>0$, we select the characteristic function for the disk centered at $(0,0)$ with radius $a>0$,
 \beq
 f_a(x) =  1 \text{ for } |x|< a \text{ and } f_a(x) = 0 \text{ for } |x| \geq a. 
 \eeq
 Taking the Fourier transform in $\mbr^2$ and using polar coordinates, we get 
 \beq
 \begin{gathered}
 \hat f_a(\xi) = \int_{\mbs^1} \int_0^a e^{-ir \omega \cdot \xi} r dr d\omega 
  = \int_{\mbs^1} \int_0^1 e^{-ir \omega \cdot a\xi} a^2 r dr d\omega = a^2 \hat f_1(a \xi).
 \end{gathered}
 \eeq
Below, we use homogeneous Sobolev spaces to compute 
 \beq
 \begin{gathered}
\|\chi_{\delta + \eps}(D)  f_a\|_{\dot H^{-1/2}}^2 = \int_{\mbr^2} |\xi|^{-1} \chi_{\delta + \eps}(\xi)^2 |\hat f_a(\xi)|^2 d\xi \\
 = \int_{\mbr^2}  |\xi|^{-1} \chi_{\delta + \eps}(\xi)^2 a^4 |\hat f_1(a \xi)|^2  d\xi 
  = \int_{\mbr^2} |\xi/a|^{-1} \chi_{\delta + \eps}(\xi/a)^2 a^2 |\hat f_1(\xi)|^2  d\xi .
\end{gathered}
 \eeq
 We see that for $a$ large, 
 \beq
 \|\chi_{\delta + \eps}(D)  f_a\|_{\dot H^{-1/2}} \simeq a^{3/2} \|\chi_{\delta + \eps}(D)  f_1\|_{\dot H^{-1/2}}.
 \eeq
Similarly, we have $ \|f_a\|_{L^2} \simeq   a \|f_1\|_{L^2}$. If $\|f_1\|_{L^2}=C_0 \|\chi_{\delta + \eps}(D) f_1\|_{H^{-1/2}}$ is true for some $C_0>0$, then for $f_a$, the condition  in \eqref{eq-dn} will be violated when $a$ is large. 
We show the test result for this example in Figure \ref{fig-disk}. As $a$ increases, we start to see some instability in the network reconstruction, and eventually, the artifacts appear around the missing angles, which agrees with our theory. 

\begin{figure}[t]
\centering
\includegraphics[scale=0.7]{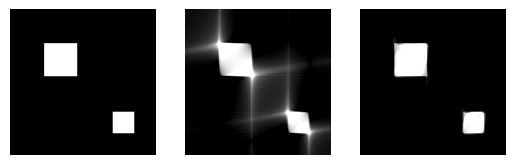}
\caption{An example of two squares (outside the training samples).Left: the original image. Middle: the FBP reconstruction with missing angle $\theta=80^{\circ}$. Right: the U-Net reconstruction, the range of pixel values for the reconstruction is $[-0.08, 1.26]$. For display purposes, all image pixel values are clipped to $[0, 1]$.}
\label{fig-square}
\end{figure}

\begin{figure}[!hbt]
\centering
\includegraphics[scale=0.7]{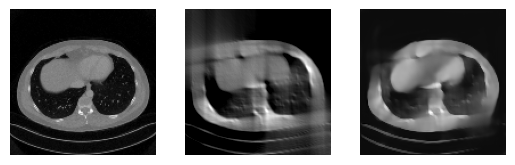}
\caption{An example from 2016 NIH-AAPM-Mayo Clinic Low Dose CT Grand Challenge.  Left: the original image. Middle: the FBP reconstruction with missing angle $\theta=80^{\circ}$. Right: the U-Net reconstruction, the range of pixel values for the reconstruction is $[-0.04, 0.88]$. For display purposes, all image pixel values are clipped to $[0, 1]$.}
\label{fig-mayo}
\end{figure}

We also test the trained network on other types of images, see an example of two squares in Figure \ref{fig-square} and an example from the 2016 NIH-AAPM-Mayo Clinic Low Dose CT Grand Challenge~\cite{mayo} in Figure~\ref{fig-mayo}, both examples show the stability can be extended beyond the visible singularities, which agree with our previous observation. It is well known that the corners in Figure \ref{fig-square} could produce strong artifacts as shown in the FBP reconstruction. Because the two squares are not too large, condition~\eqref{eq-dn} is likely to hold via a similar Fourier analysis as for $f_a$.  
 
  %========================%
\subsection{Further discussions} \label{sec:dis}
\begin{figure}[!htb]
    \centering
    \includegraphics[width=0.8\linewidth]{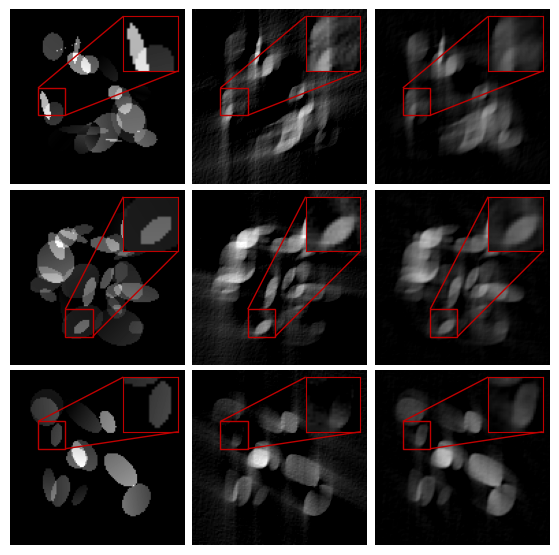}
    \caption{Three numerical reconstructions using FBP and U-Net without nonlinear activation functions. These examples are the same as Figure~\ref{fig-test}. Each image's top-right shows a zoomed-in detail. From top to bottom: missing angle $\theta=100^{\circ}, 80^{\circ}, 60^{\circ}$. Left: the original image. Middle: the FBP reconstruction. Right: the U-Net reconstruction without nonlinear activations, from top to bottom, the ranges of pixel values of the reconstructions are $[-0.288, 0.671]$, $[-0.124, 0.734]$, $[-0.208, 0.767]$. For display purposes, all image pixel values are clipped to $[0, 1]$. }
    \label{fig:unet linear}
\end{figure}

 The examples presented in Sections \ref{subsec-net} and \ref{subsec-rob} indicate that the trained neural network pre-composed with the FBP could be a good approximation of $R_\delta^{-1}$. However, we do not claim that every trained network shares the same property. To prevent overfitting, we stop the training after a small number of iterations, see Figure~\ref{fig-train_loss}. Although training the network for a longer time may keep training errors decreasing, the test errors on the test set and general images in Section \ref{subsec-rob} usually worsen.  We also emphasize that the choice of activation function is essential in the U-Net. First, the activation has to be an affine mapping. Without the nonlinearity in the activations, the effect of the convolution layers can be viewed as a single convolution kernel $\mathcal{K}$. Assuming the noiseless data that $F_j = R_{\delta} f_j$,  $j=1,\cdots, n$, then using Plancherel's theorem, the loss function becomes 
\begin{equation}
\begin{aligned}
 &\frac{1}{n}\sum_{j=1}^{n} \left\|\int_{\mathbb{R}^2} \mathcal{K}(x - y) \mathcal{A}^{\dagger} F_j(y) dy - f_j(x) \right\|_{L^2}^2 \\
 =& \frac{1}{n} \sum_{j=1}^{n} \left\|\widehat{\mathcal{K}}(\zeta) \widehat{f_j}(\zeta)\chi_{\delta}(\zeta) - \widehat{f_j}(\zeta) \right\|_{L^2}^2.
\end{aligned}
\end{equation}
It implies that the minimization of the loss function cannot recover any information outside the cone $\mathcal{K}_{\delta}$, see Figure~\ref{fig:unet linear} for the reconstructions for the same examples as Figure~\ref{fig-test} but using U-Net without nonlinear activation functions. 

To demonstrate the potential effectiveness of the activation function, we consider a simplified shallow U-Net that skips all downsampling and upsampling blocks (only consists of the bottleneck and output layer). For analysis purposes, we minimized the loss function by gradient flow.
\begin{equation}
\begin{aligned}
   \mathcal{N}_{\boldsymbol{\theta}} g:=  \sum_{i=1}^n a_i\sigma\left(\sum_{j=1}^{n}  K_{ij} \ast g\right) = \sum_{i=1}^n\tilde{a}_i\sigma\left(\frac{1}{\sqrt{n}}\sum_{j=1}^{n}  K_{ij} \ast g\right), 
\end{aligned}
\end{equation}
where $\tilde{a}_i:= \sqrt{n} a_i$, $\sigma(x) =  \max(0, x) +\alpha\min(0, x)$ with typical value $\alpha = 10^{-2}$, and each $K_{ij}$ is a convolutional kernel function. The loss function becomes 
\begin{equation}
    \mathcal{L} = \frac{1}{|\mathcal{D}|} \sum_{f\in \mathcal{D}} \int_{\mathcal{X}}| \mathcal{N}_{\boldsymbol{\theta}  }\mathcal{A}^{\dagger} R_{\delta} f(x) - f(x)|^2  dx, 
\end{equation}
where $\mathcal{X} =[0, 1]^2$ denotes the support of the data set $\mathcal{D}$. Initially, the \emph{band-limited} kernel functions $\{K_{ij}\}_{1\le i, j\le n}$ are generated from the random field 
$$K(x, \omega) = \sigma_K \sum_{|\zeta|\le B} Z_{\zeta}(\omega) e^{i2\pi \zeta\cdot x}, $$
where $B$ denotes a sufficiently wide band limit, $\sigma_K$ is a configuration constant, and $Z_{\zeta}(\omega) = \overline{Z_{-\zeta}(\omega)}$ are i.i.d standard normal variables (real and imaginary parts) on the half-plane, except $|\zeta| = 0$. In the regime of the neural tangent kernel (NTK)~\cite{ntk, bias} that $n\to\infty$, the kernel functions $K_{ij}$ can be regarded as frozen and the only trainable parameters become $\{\tilde{a}_i\}_{i=1,\cdots, n}$. For any $f\in \mathcal{D}$, 
 \begin{equation}
 \begin{aligned}
     K(x, \omega)\ast \mathcal{A}^{\dagger} R_{\delta} f& = \sigma_K \sum_{ |\zeta|\le B} Z_{\zeta}(\omega) e^{i2\pi\zeta\cdot  x} \ast \mathcal{A}^{\dagger} R_{\delta} f  \\&= \sigma_K \sum_{ |\zeta|\le B }Z_{\zeta}(\omega)\widehat{f}(\zeta)\chi_{\delta}(\zeta) e^{i2\pi \zeta \cdot x}, 
 \end{aligned}
\end{equation}
which defines another random field, denoted by $K_{\delta}(x, f, \omega)$, with correlation function 
\begin{equation}
    \rho_{\delta}(x-y, f, h) = \frac{1}{ \|f\|_{\delta}  \|h\|_{\delta}} \sum_{|\zeta|\le B, \zeta\in \mathcal{K}_{\delta}} \widehat{f}(\zeta) \overline{\widehat{h}(\zeta)} e^{i2\pi \zeta\cdot(x - y)},
\end{equation}
where $
    \|f\|_{\delta}:=  \sqrt{\sum_{\zeta\in \mathbb{Z}^2\cap \mathcal{K}_{\delta}} |\widehat{f}(\zeta)|^2}$. 
Without the activation function, this correlation has a large kernel space containing Fourier modes $e^{i2\pi \zeta\cdot x}$ that $\zeta\notin \mathcal{K}_{\delta}$.
After applying the nonlinear activation function $\sigma$, we denote $\widetilde{\rho}_{\delta}$ as the new correlation function for the transformed field $\sigma(K_{\delta}(x, f, \omega))$, which satisfies
\begin{equation}
   \widetilde{\rho}_{\delta} := \rho_{\delta} + \frac{(1-\alpha)^2}{\pi (1 + \alpha^2)} (\sqrt{1-\rho_{\delta}^2} - \rho_{\delta}\cos^{-1}(\rho_{\delta})).
\end{equation}
Then, the gradient flow for the loss function $\mathcal{L}$ becomes 
\begin{equation}
    \frac{d\mathcal{L}}{dt} = -\frac{2\sigma_K^2}{|\mathcal{D}|}\sum_{f, h\in \mathcal{D}} \int_{\mathcal{X}\times \mathcal{X}} T_{\boldsymbol{\theta}}f(x)\widetilde{\rho}_{\delta} (x - y, f, h)   T_{\boldsymbol{\theta}}h(y) dx dy,
\end{equation}
where $T_{\boldsymbol{\theta}} f = \|f\|_{\delta}(N_{\boldsymbol{\theta}} f - f)$. Moreover, if we decompose the correlation function into Fourier basis
\begin{equation}
    \widetilde{\rho}_{\delta}(x -y, f, h) = \sum_{\zeta\in \mathbb{Z}^2} C_{\zeta}(f, h) e^{i 2\pi \zeta \cdot (x - y)}, 
\end{equation}
then if each $|\mathcal{D}|\times |\mathcal{D}|$ matrix $C_{\zeta}(f, h)$, $f, h\in \mathcal{D}$ is strictly positive definite, then the loss function $\mathcal{L}$ will eventually converge to zero. 

As a special example, if the data set $\mathcal{D}$ satisfies that $\rho_{\delta}(s, f, h) \equiv 0$ for any pair $f\neq  h\in \mathcal{D}$, then we only have to focus on the diagonal part $\rho_{\delta}(s, f, f)$.  When $|s|\ll 1$, we can expand $\rho_{\delta}(s, f, f) = 1 - \frac{A}{2}s^2 + \mathcal{O}(s^4)$ for a constant $A > 0$ depending on $f$, which means 
 \begin{equation}
     \widetilde{\rho}_{\delta}(s, f, f) = {\rho}_{\delta}(s, f, f) + \frac{(1-\alpha)^2}{\pi (1 + \alpha^2)} \frac{A^{3/2}}{3} |s| s^2 + \mathcal{O}(s^4).
 \end{equation}
The new kernel function finds a mild non-smoothness at $s=0$, which implies that if $|\zeta|$ is sufficiently large, the Fourier modes $e^{i2\pi \zeta\cdot x}$ cannot stay inside the kernel space. Hence, the nonlinear activation enables the network to fit some information beyond the cone $\mathcal{K}_{\delta}$ under the NTK setting. Nevertheless, in this regime, the network often has severe overfitting, hence the early stopping strategy is necessary. The discussion of the training mechanism in a more general setting is beyond the scope of this work.

 \section{Conclusion}
 \label{sec-conclude}
 In this note, we provided a mathematical explanation of the robustness of data-driven approaches based on deep learning for limited-angle tomography. We showed that prior conditions on the training data set can lead to a stable reconstruction of some invisible singularities. We also tested the theory using a simple U-Net with a pre-composed FBP. For data sets $\mathcal{D}\subset \mathcal{D}_{N,\eps}$ (see~\eqref{eq-dn}) with a relatively small $N$ (e.g. the public data set from~\cite{BGR}), the numerical experiments demonstrate stable reconstruction beyond the visible angle. For data sets with growing $N$ (e.g. Figure~\ref{fig-disk}), the reconstructed images suffer severely from the artifacts. These numerical experiments agree with our analysis in Section 3.

\section*{Acknowledgment}
The authors would like to thank the anonymous referees who provided useful and detailed comments. YW is partly supported by NSF under grant DMS-2205266. YZ is partly supported by NSF under grant DMS-2309530.

%%===============================REFERENCE==========================================%

\end{document}